\documentclass{elsart}

\usepackage{graphics}
\usepackage{graphicx}
\usepackage{epsfig}

\usepackage{amssymb}
\usepackage{amsmath}

\newcommand{\di}{\mathrm{div}}
\newcommand{\p}{\partial}

\newcommand{\norm}[1]{\left\Vert#1\right\Vert}

\newcommand{\supp}{\mathrm{supp}}
\makeatletter

\newcommand{\Rmnum}[1]{\expandafter\@slowromancap\romannumeral#1@}
\makeatother

\begin{document}

\begin{frontmatter}



\title{Transonic Potential Flows
in A Convergent--Divergent Approximate Nozzle}


\author[Yuan]{Hairong Yuan\thanksref{yuan}}
\author[He]{Yue He\thanksref{he}}

\address[Yuan]{
Department of Mathematics, East China Normal University, Shanghai
200241, China} \ead{hryuan@math.ecnu.edu.cn \&
hairongyuan0110@gmail.com}

\address[He]{
Department of Mathematics, School of Mathematics and Computing,
Nanjing Normal University, Nanjing 210097, China}
\ead{heyueyn@163.com \& heyue@njnu.edu.cn}

\thanks[yuan]{Supported by  China Postdoctoral Science
Foundation (20070410170), and DMS-0720925 ``U.S.-China CMR:
Multidimensional Problems in Nonlinear Conservation Laws and Related
Applied Partial Differential Equations" by U.S. NSF and Chinese
NSF.}

\thanks[he]{Supported by NNSF grant of P. R. China:
No.\,10571087, and Natural Science Foundation of Jiangsu Education
Commission of P. R. China: No.\,06KJB110056.}

\begin{abstract}
In this paper we prove existence, uniqueness and regularity of
certain perturbed (subsonic--supersonic) transonic potential flows
in a two-dimensional Riemannian manifold with
``convergent--divergent" metric, which is an approximate model of
the de Laval nozzle in aerodynamics. The result indicates that
transonic flows obtained by quasi-one-dimensional flow model in
fluid dynamics are stable with respect to the perturbation of the
velocity potential function at the entry (i.e., tangential velocity
along the entry) of the nozzle. The proof is based upon linear
theory of elliptic--hyperbolic mixed type equations in physical
space and a nonlinear iteration method.

\end{abstract}

\end{frontmatter}


\section{Introduction}

Understanding flow patterns in a convergent--divergent nozzle (the
so called de Laval nozzle in engineering) is a prominent issue in
aerodynamics and partial differential equations  due to their
numerous applications in practice, and closely connection with many
difficult mathematical problems, such as mixed type equations and
free boundary problems \cite{CF,K1,K2,Mo2}. Since these flow
patterns are genuinely nonlinear, various physically significant
special solutions of the corresponding mathematical problems play an
important role in the theoretical analysis. For example, for the
nearly spherical symmetric transonic shocks and transonic shocks in
a slowly varying nozzle, there are works of Chen et.al.
\cite{CF1,CY,LY,XY,Yu2} on potential flow equation and complete
Euler system based  upon in essence two classes of  special
solutions \cite{Yu3}. There are also many progresses in the analysis
of subsonic nozzle flows, see, for instance, \cite{L,XX} and
reference therein. However, since  no simple and physical special
transonic--flow solution is available, presently the study of
subsonic--supersonic transonic flow  mainly utilized the methods of
compensated compactness (see \cite{CDSW,CSW,Mo2,XX} and references
therein).

In \cite{Yu1}, motivated by a significant work of Sibners \cite{SS},
Yuan constructed various interesting special solutions in a
two--dimensional Riemannian manifold with ``convergent--divergent"
metric, which may be regarded as an approximation of the de Laval
nozzle. This paper also studied several boundary value problems of
subsonic flows in such a manifold. In the present paper we will
further investigate the subsonic--supersonic flow via the potential
flow equation. We show that the special subsonic--supersonic
transonic flows are stable with respect to the small perturbation of
the velocity potential function (i.e., the tangential velocity) at
the entry (see Theorem \ref{thm301}).

The potential flow equation is a second order equation of
elliptic--hyperbolic mixed type for transonic flow. For such
equations, presently one of the main tool is the theory of positive
symmetric systems and techniques of energy estimates, see, for
instance, \cite{HH,K1,K2} and references therein. In this paper we
employ the theory developed in \cite{K2} to show the solvability of
linear problem, and then a nonlinear iteration argument to solve the
nonlinear problem.

We remark that in recent years there are many breakthroughs on
partial differential equations of mixed type and degenerate elliptic
type arising in differential geometry and physics, see, for example,
\cite{CF2,HH,He1,He2,Kim1,Kim2,K2,Mo2}. For earlier developments in
this field, one may also consult \cite{Mo2,S1,S2} and references
therein.

The rest of the paper is organized  as follows. In Section 2 we
formulate the problem, and study the properties of the coefficients
of the potential flow equations in the manifold. In Section 3 we
solve the linear problem, and finally in Section 4 we state the main
result, Theorem \ref{thm301} and  prove it.

\section{Formulation of the problem}

Let $\mathbf{S}^1$ be the standard unit circle in $\mathbf{R}^2$,
and $\mathcal{M}$ be the Riemannian manifold
$\{(x^1,x^2)\in[-1,1]\times\mathbf{S}^1\}$ with a metric
$G=g_{ij}dx^i\otimes dx^j=dx^1\otimes dx^1+n(x^1)^2dx^2\otimes
dx^2.$ Here $n(t)$ is a positive smooth function on $[-1,1]$
satisfies: (1) $n''(t)>0;$ (2) $n'(t)<0$ for $t\in(-1,0)$, $n'(t)>0$
for $t\in(0,1).$ Such a manifold $\mathcal{M}$ may be regarded as an
approximation of a two--dimensional convergent--divergent nozzle,
with $\mathcal{M}^\pm=\mathcal{M}\cap\{x^1\gtrless0\}$ respectively
the divergent and convergent part. We also call
$\Sigma^{k}=\{k\}\times \mathbf{S}^1,\ k=-1,0,1 $ respectively the
entrance, throat and exit of $\mathcal{M}.$ Obviously
$\p\mathcal{M}=\Sigma^{-1}\cup\Sigma^1.$

Let $p,\rho$ be functions in $\mathcal{M}$ represent respectively
the pressure and density of  gas flow in $\mathcal{M}$, and $v$ be a
vector field in $\mathcal{M}$ represent the velocity of the flow. We
consider polytropic and isentropic gas flows, then
$p=\kappa\rho^\gamma$ with $\kappa>0,\ \gamma>1$ two constants, and
the speed of sound is $c=\sqrt{\kappa\gamma\rho^{\gamma-1}}.$ Let
$\bar{v}$ be the 1-form corresponding to $v$ under the metric $G$.
The flow is irrotational if $\bar{v}$ is exact; That is, there
exists a function $\varphi$ in $\mathcal{M}$ such that
$\bar{v}=d\varphi.$ Substituting this in the equation of
conservation of mass $\di (\rho v)=-d^*(\rho\bar{v})=0,$ where $\di$
and $d^*$ are respectively the divergence operator and
codifferential operator in $\mathcal{M},$ then by the formula
$d^*(\rho \bar{v})=\rho d^*\bar{v}-\langle d\rho, \bar{v}\rangle,$
with $\langle\cdot,\cdot\rangle$ the inner product of forms in
$\mathcal{M}$,  we have
\begin{eqnarray}\label{101}
\rho \Delta\varphi=\langle d\rho,d\varphi\rangle,
\end{eqnarray}
where $\Delta=dd^*+d^*d$ is the Hodge Laplacian of forms. (Note that
$d^*\varphi=0.$) By the Bernoulli's law which represents
conservation of energy:
\begin{eqnarray}\label{102}
\frac 12\langle
d\varphi,d\varphi\rangle+\frac{\kappa\gamma}{\gamma-1}\rho^{\gamma-1}=c_0,
\end{eqnarray}
where $c_0$ is a positive constant, $\rho$ in \eqref{101} can be
expressed in terms of $d\varphi.$ So we may write \eqref{101} as a
second order equation of $\varphi$.

Indeed, in the $(x^1,x^2)$ coordinates, we have
\begin{eqnarray*}
\Delta\varphi&=&-\frac{1}{\sqrt{g}}\p_i(\sqrt{g}g^{ij}\p_j\varphi)\\
&=&-\frac{1}{n(x^1)}\left(\p_1(n(x^1)\p_1\varphi)+\frac{1}{n(x^1)}\p_{22}\varphi\right).
\end{eqnarray*}
Here $\sqrt{g}=\sqrt{\det(g_{ij})},$\  $(g^{ij})$ is the inverse of
$(g_{ij}),$ and $\p_i=\p_{x^i}$, \ $\p_{ij}=\p_i\p_j.$ By
differentiating \eqref{102} we have
\begin{eqnarray*}
\frac{c^2}{\rho}d\rho=-\big(\frac 12
\p_i\varphi\p_j\varphi\p_kg^{ij}+\p_{ik}\varphi\p_j\varphi
g^{ij}\big)dx^k.
\end{eqnarray*}
Then by a straightforward calculation we obtain
\begin{eqnarray}\label{103}
&&n(x^1)^2(c^2-(\p_1\varphi)^2)\p_{11}\varphi-2\p_1\varphi\p_2\varphi\p_{12}\varphi
+\left(c^2-\frac{1}{n(x^1)^2}(\p_2\varphi)^2\right)\p_{22}\varphi\nonumber\\
&&\qquad\qquad
+n(x^1)n'(x^1)\left(c^2+\frac{1}{n(x^1)^2}(\p_2\varphi)^2\right)\p_1\varphi=0.
\end{eqnarray}
Direct computation yields that this equation is of elliptic type if
the flow is subsonic
($c^2>(\p_1\varphi)^2+(\p_2\varphi)^2/n(x^1)^2),$ and is of
hyperbolic type if the flow is supersonic
($c^2<(\p_1\varphi)^2+(\p_2\varphi)^2/n(x^1)^2).$

If the flow depends only on $x^1,$ then \eqref{103} indicates that
$\varphi_b=\varphi_b(x^1)$ satisfies the equation
\begin{eqnarray}\label{104}
&&n(x^1)(c_b^2-(\p_1\varphi_b)^2)\p_{11}\varphi_b
+n'(x^1)c_b^2\p_1\varphi_b=0,
\end{eqnarray}
where $c_b$ is the sonic speed corresponding to $\varphi_b.$ It can
be shown that there are special flows $\varphi_b\in
C^5(\mathcal{M})$ which are subsonic in $\mathcal{M}^-$ and
supersonic in $\mathcal{M}^+$, and $\p_1\varphi_b>0,\
\p_{11}\varphi_b>0$ in $\mathcal{M}$ (see \cite{Yu1}). We call such
flows {\it background solutions.} The aim of this paper is to study
stability of certain background solutions under perturbations of
$\varphi$ on the entrance $\Sigma^{-1}$ of $\mathcal{M}.$

Let $\hat{\varphi}=\varphi-\varphi_b.$ By subtracting \eqref{104}
from \eqref{103}, we have
\begin{eqnarray}\label{105}
k(D\varphi)\p_{11}\hat{\varphi} +b(D\varphi)\p_{12}\hat{\varphi}
+\p_{22}\hat{\varphi}-\alpha(x^1)\p_1\hat{\varphi}=f(D{\varphi}),
\end{eqnarray}
where
\begin{eqnarray}
k(D\varphi)&:=&\frac{n(x^1)^2(c^2-(\p_1\varphi)^2)}{c^2-\frac{1}{n(x^1)^2}(\p_2\varphi)^2},\label{106}\\
b(D\varphi)&:=&-\frac{2\p_1\varphi\p_2\varphi}{c^2-\frac{1}{n(x^1)^2}(\p_2\varphi)^2},\label{107}\\
\alpha(x^1)&:=&\frac{n(x^1)^2\p_{11}\varphi_b}{\p_1\varphi_b}\cdot\frac{c_b^4+c_b^2(\p_1\varphi_b)^2
+(\gamma-1)(\p_1\varphi_b)^4}{c_b^4},\label{109}
\end{eqnarray}
\begin{eqnarray}
f(D{\varphi})&:=&\left\{\frac{\gamma-1}{2}\Big(\p_{11}\varphi_b+\frac{n'(x^1)}{n(x^1)}
\p_1\varphi_b\Big)\left((\p_1\hat{\varphi})^2
+\frac{1}{n(x^1)^2}(\p_2\hat{\varphi})^2\right)\right.\nonumber\\
&&\left.+\p_{11}\varphi_b(\p_1\hat{\varphi})^2-\frac{n'(x^1)}{n(x^1)}(c^2-c_b^2)\p_1\hat{\varphi}
-\frac{n'(x^1)}{n(x^1)^3}\p_1\varphi(\p_2\hat{\varphi})^2\right\}\nonumber\\
&&\qquad\cdot\frac{n(x^1)^2}{c^2-\frac{1}{n(x^1)^2}(\p_2\varphi)^2}\\
&&+\left[\frac{n(x^1)^2\p_{11}\varphi_b}{c_b^2\p_1\varphi_b}\cdot\frac{c_b^4+c_b^2(\p_1\varphi_b)^2
+(\gamma-1)(\p_1\varphi_b)^4}{c^2-\frac{1}{n(x^1)^2}(\p_2\varphi)^2}-\alpha(x^1)\right]\p_1\hat{\varphi}
.\label{110}\nonumber
\end{eqnarray}
We will investigate the following problem:
\begin{eqnarray}
&\text{Eq.}\ \  \eqref{105} & \text{in}\ \  \mathcal{M}\ \
\text{with}\ \
\varphi=\hat{\varphi}+\varphi_b,\label{111}\\
&\hat{\varphi}=g(x^2)& \text{on}\ \  \Sigma^{-1}\ \  \text{with} \ \
\norm{g}_{H^5(\mathbf{S}^1)}\ \  \text{small.}\label{112}
\end{eqnarray}
This is a Dirichlet problem of an elliptic--hyperbolic mixed
equation.

For background solution $\varphi_b,$ let
$\tau=(\p_1\varphi_b/c_b)^2$ be the square of Mach number. By
Bernoulli's law we may compute
\begin{eqnarray}\label{113}
\p_1(k(D\varphi_b))=-\frac{nn'}{\tau-1}[(\gamma+1)\tau^2-2\tau+2]<0
\end{eqnarray}
in $\mathcal{M}$. Here we also note the fact that
$\p_{11}\varphi_b=(\p_1\varphi_b)n'/(n\cdot(\tau-1)),$ and
$n'/(\tau-1)>0$ in $\mathcal{M}$, especially at the throat
$\Sigma^0$ we have $n'/(\tau-1)=\sqrt{nn^{''}/(\gamma+1)}$ by
L'Hospital's rule in calculus. We also may write
\begin{eqnarray*}
\alpha(x^1)=\frac{nn'[1+\tau+(\gamma-1)\tau^2]}{\tau-1}>0,
\end{eqnarray*}
so there holds
\begin{eqnarray}
2\alpha-l\p_1(k(D\varphi_b))>\delta_1>0\label{114}
\end{eqnarray}
in $\mathcal{M}$ for a fixed number $\delta_1$ and all positive
number $l.$  In addition,
\begin{eqnarray}\label{115}
2\alpha+\p_1(k(D\varphi_b))=\frac{nn'\tau}{\tau-1}[(\gamma-3)\tau+4]>\delta_2>0
\end{eqnarray}
in $\mathcal{M}.$ We remark that the constants $\delta_1,\ \delta_2$
depend only on the specific background solution $\varphi_b$.

\section{Solvability of linear problem}

Let $k,\ b,\ a,\ \alpha,\ f$ be functions in
$M=[-1,1]\times\mathbf{S}^1$ with flat metric. In this section we
investigate the solvability of the following linear problem:
\begin{eqnarray}
&Lu:=k\p_{11}u+b\p_{12}u +a\p_{22}u-\alpha(x^1)\p_1u=f &\text{in}\ \
{M},\label{201}\\
&u=0 &\text{on}\ \ \Sigma^{-1}.\label{202}
\end{eqnarray}
We use $H^j$ to denote the Sobolev space $W^{j,2}({M}),$ and the
corresponding norm is written as $\norm{\cdot}_j$. Note that
$C^{j}({M})$ denotes the usual space of $j$ times continuous
differential functions in ${M}$.

\begin{thm}\label{thm201}
Let $k>0$ on $\Sigma^{-1}$ and $k<0$ on $\Sigma^1,$ and $k, b, a,
\alpha\in C^4{({M})}.$ Suppose $f\in H^s({M}), 0\le s\le3$ and there
is a positive constant $\delta$ such that
\begin{eqnarray}
&a\ge\delta>0,\quad \p_1a\le-\delta<0 &\text{in}\ \ {M},\label{h0}\\
&2\alpha-(2p-1)\p_{1}k\ge\delta>0 &\text{in}\ \ {M}\ \ \text{for}\ \
p=0,1,\cdots, s.\label{203}
\end{eqnarray}
Then there is a $\nu>0$ such that if
\begin{eqnarray}\label{2005}
\norm{b}_{3}\le \nu,\ \  \norm{\p_2a}_{C^3}\le \nu \qquad \text{in}\
\ {M},
\end{eqnarray} then  there  exists uniquely one solution $u\in H^{1}$ to
problem \eqref{201} and \eqref{202} and there holds the estimate
\begin{eqnarray}\label{204}
\norm{u}_{s+1}\le C\norm{f}_s.
\end{eqnarray}
\end{thm}

To prove this, we follow the celebrated ideas presented in
\cite{K1,K2}. That is, one first solves a boundary value problem of
a mixed type equation which is elliptic at both the entry and exit.
Then the above theorem can be demonstrated by an extension
technique.

\begin{thm}\label{thm302}
Under the assumptions of Theorem \ref{thm201}, but supposing $k>0$
on both $\Sigma^{-1}$ and $\Sigma^1$, then the following problem
\begin{eqnarray}
&Lu=f &\text{in} \ \ {M},\label{h1}\\
&u=0 &\text{on} \ \ \Sigma^{-1},\label{h2}\\
&\p_1u=0 &\text{on} \ \ \Sigma^{1}\label{h3}
\end{eqnarray}
has uniquely one solution $u\in H^1$ and it also satisfies
\eqref{204}.
\end{thm}

\noindent{\it Proof of Theorem \ref{thm201}.} {\it Step 1.\
Uniqueness.} Let $f\equiv0$. Multiplying  to \eqref{201} by $\p_1u$
and integrating the expression in $[-1,1]\times\mathbf{S}^1$, Note
that $\p_2u=0$ on $\Sigma^{-1}$, we have
\begin{eqnarray}
&&\int_{M}\big[(\alpha+\frac12 \p_1k+\frac12
\p_2b)(\p_1u)^2-\frac12\p_1a(\p_2u)^2+\p_2a\p_1u\p_2u\big]\, dx^1dx^2\nonumber\\
&=&\int_{\Sigma^1}\big[\frac{k}2(\p_1u)^2-\frac
a2(\p_2u)^2\big]\,dx^2-\int_{\Sigma^{-1}}\big[\frac{k}2(\p_1u)^2-\frac
a2(\p_2u)^2\big]\,dx^2\nonumber\\
&\le&0 \label{3007}
\end{eqnarray}
via the integration by parts and divergence theorem. Since
$\p_1a\le-\delta<0$ in $\mathcal{M}$, we infer that, by choosing
$\nu=\delta/4$,
\begin{eqnarray}
&&(2\alpha+\p_1k+\p_2b)(\p_1u)^2-\p_1a(\p_2u)^2+2\p_2a\p_1u\p_2u\nonumber\\
&\ge&(\delta-2\nu)(\p_1u)^2+(\delta-\nu)(\p_2u)^2
\ge\frac{\delta}{2}((\p_1u)^2+(\p_2u)^2).\nonumber
\end{eqnarray}
Hence we have
\begin{eqnarray}\label{h4}
(\delta/2)\int_{{M}}[(\p_1u)^2+(\p_2u)^2]\,dx^1dx^2\le0
\end{eqnarray}
by \eqref{3007} and therefore $u\equiv0$ in ${M}$ due to
\eqref{202}.

{\it Step 2.\ Existence.} Let ${M}^*=[-1,2]\times\mathbf{S}^1.$ We
may extend $k,b,a,\alpha,f$ to ${M}^*$ such that  they still satisfy
\eqref{203}--\eqref{2005} and other assumptions in Theorem
\ref{thm302}, especially $k>0$ on the new exit
$\Sigma^2=\{2\}\times\mathbf{S}^1,$ and $\norm{f}_{H^s(M^*)}\le
C\norm{f}_{H^s(M)}.$ Denoting the obtained extended operator in
${M}^*$ as $L$ again, we consider the problem \eqref{h1}\eqref{h2}
together with boundary condition $\p_1u=0$ on $\Sigma^2.$ By Theorem
\ref{thm302}, there is a unique solution $u^*$ and
$\norm{u^*}_{s+1}\le C\norm{f}_{H^s(M^*)}.$ Obviously
$u:=u^*|_{{M}}$ is also a solution to problem \eqref{201} and
\eqref{202}. This finishes the proof.

\smallskip
The proof of Theorem \ref{thm302} also follows in a similar way of
\cite{K2} (Theorem 1.1, pp.9--18), but needs some modifications to
deal with the mixed derivative term $b\p_{12}u$, the non-divergence
term $a\p_{22}u$, and no lateral boundary in our case. For
completeness and convenience of the readers, we sketch out the
proofs. Some of the details in the proof are important in the
analysis of the nonlinear problems.

\smallskip
\noindent{\it Proof of Theorem \ref{thm302}.} {\it Step 1.\
Uniqueness.} This may be proved by a similar method as in deriving
\eqref{h4}.

{\it Step 2.\ Approximate problem.} To show existence of a $H^1$
weak solution to problem \eqref{h1}--\eqref{h3}, as in \cite{K2}, we
consider the following singular perturbation problem:
\begin{eqnarray}
L_\epsilon u:=Lu+\epsilon\p_{111}u=f_\epsilon && 1>\epsilon>0\ \
\text{in}\ \
{M},\label{3009}\\
u(-1,x^2)=\p_1u(-1,x^2)=0,&&\label{3010}\\
\p_1u(1,x^2)=0,&&\label{3011}
\end{eqnarray}
where $f_\epsilon\in C^1$ and $f_\epsilon\rightarrow f$ in $L^2$ as
$\epsilon\rightarrow0.$

{\it Step 2.1.\ Galerkin's method.} To show existence of a solution
$u^\epsilon$ to problem \eqref{3009}--\eqref{3011}, we employ the
Galerkin's method of finite dimensional approximation. Let
$\{Y_i(x^2)\}\ (i=1,2,\cdots)$ be a complete system in
$H^2(\mathbf{S}^1)$ and orthogonal in $L^2(\mathbf{S}^1)$. We may
also assume that each $Y_i(x^2)$ is smooth. Set
\begin{eqnarray*}
u^{N,\epsilon}(x^1,x^2)=\sum_{i=1}^N X_i^{N,\epsilon}(x^1)Y_i(x^2),
\ \ N=1,2,\cdots.
\end{eqnarray*}
The functions $X_i^{N,\epsilon}\ (i=1,\cdots,N)$ are to be
determined by a boundary value problem of a system of  third order
ODEs:
\begin{eqnarray}
\int_{\mathbf{S}^1}[L_\epsilon
u^{N,\epsilon}-f_\epsilon]Y_j(x^2)\,dx^2=0,&&\
\ j=1,2,\cdots,N,\label{313}\\
X_i^{N,\epsilon}(-1)=(X_i^{N,\epsilon})'(-1)=0,&&\label{314}\\
(X_i^{N,\epsilon})'(1)=0.&&\label{315}
\end{eqnarray}
Equation \eqref{313} can be written explicitly as
\begin{eqnarray}
&&\epsilon (\int_{\mathbf{S}^1}dx^2)(X_j^{N,\epsilon})'''
+\sum_{i=1}^N\left[\left(\int_{\mathbf{S}^1}kY_iY_j\,dx^2\right)
(X_i^{N,\epsilon})''\right.\nonumber\\
&&\quad\left. \right.+\left(\int_{\mathbf{S}^1}(bY_i'Y_j-\alpha
Y_iY_j)\,dx^2\right)(X_i^{N,\epsilon})'
\nonumber\\
&&\qquad \ \
\left.+\left(\int_{\mathbf{S}^1}aY_i''Y_j\,dx^2\right)(X_i^{N,\epsilon})\right]
=\int_{\mathbf{S}^1}f_\epsilon Y_j(x^2)\,dx^2.\label{316}
\end{eqnarray}

{\it Step 2.1.1.\ Uniqueness.} Now we show the solution to problem
\eqref{313}--\eqref{315} is unique. Indeed, multiplying to
\eqref{313} by $(X_j^{N,\epsilon})'$,  summing up for $j$ from $1$
to $N$ and integrating with respect to $x^1$ on $[-1,1]$, we have
\begin{eqnarray*}
-\int_{M}[L_\epsilon  u^{N,\epsilon}
\p_1u^{N,\epsilon}]\,dx^1dx^2=-\int_{M}f_\epsilon\p_1u^{N,\epsilon}\,dx^1dx^2.
\end{eqnarray*}
Writing $u^{N,\epsilon}$ simply as $w$, then integrating by parts
and using $\p_1w(-1,x^2)=\p_1w(1,x^2)=0$, we obtain that
\begin{eqnarray}
&&-2\int_{M}[L_\epsilon  w \p_1w ]\,dx^1dx^2
=\int_{M}[(2\alpha+\p_1k+\p_2b)(\p_1w)^2-(\p_1a)(\p_2w)^2
\nonumber\\
&&\quad+2\p_2a\p_1u\p_2u]\,dx^1dx^2+2\epsilon\int_{{M}}(\p_{11}w)^2\,dx^1dx^2
+\int_{\Sigma^{-1}}[a(\p_2u)^2]\,dx^2.\nonumber
\end{eqnarray}
Therefore by H\"older inequality and Young inequality, due to
\eqref{2005}, we have
\begin{eqnarray}\label{319}
\epsilon\norm{\p_{11}w}^2_0+\norm{D w}_0^2\le
C\norm{f_\epsilon}_0^2\le C'\norm{f}_0^2.
\end{eqnarray}
Since $w(-1,x^2)=0$, it follows that \begin{eqnarray*}
w(x^1,x^2)=w(x^1,x^2)-w(-1,x^2)=\int_{-1}^{x^1}\p_1 w(t,x^2)dt.
\end{eqnarray*}
By H\"older inequality,
\begin{eqnarray*}
|w(x^1,x^2)|^2\leq \int_{-1}^{x^1}
1^2dt\,\int_{-1}^{x^1}|\p_1w(t,x^2)|^2dt \leq 2\,\int_{-1}^1
|\p_1w(x^1,x^2)|^2dx^1.
\end{eqnarray*}
Hence, we have
\begin{eqnarray}
\|w\|_0^2&=&\int_{M}|w(x^1,x^2)|^2\,dx^1dx^2\nonumber\\
&\leq& 4\,\int_{M}
|\p_1w(x^1,x^2)|^2dx^1dx^2=4\|\p_1w\|_0^2.\label{8-1-18-12-45}
\end{eqnarray}
Now it easily follows from \eqref{319} and \eqref{8-1-18-12-45} that
\begin{eqnarray}\label{320}
\epsilon\norm{\p_{11}w}^2_0+\norm{w}_1^2\le C\norm{f}_0^2
\end{eqnarray}
for a constant $C$ independent of $N$ and $\epsilon.$ This proves
uniqueness of solution to problem \eqref{313}--\eqref{315}. Note the
above estimate also holds for the solution $u^\epsilon$ of problem
\eqref{3009}--\eqref{3011}.

{\it Step 2.1.2.\ Existence and regularity.} Now by Fredholm
alternative of boundary value problems of ODEs, we can infer that
problem \eqref{313}--\eqref{315} has a solution $u^{N,\epsilon}$
which satisfies \eqref{320}. Indeed, we may write
\eqref{313}--\eqref{315} equivalently as a boundary value problem of
a first order system of ODEs with $3n$ unknowns
\begin{eqnarray}
&&\frac{d\mathcal{X}}{dx^1}=\mathcal{A}\mathcal{X}+\mathcal{F},\\
&&\mathcal{B}\mathcal{X}=0.
\end{eqnarray}
Here $\mathcal{X}, \mathcal{F}\in\mathbf{R}^{3n}$ are column
vectors, and $\mathcal{A},\mathcal{B}$ are $3n\times3n$ matrices.
Let $\mathcal{X}_1,\cdots,\mathcal{X}_{3n}$ be  a set of $3n$
linearly independent solutions of the homogeneous equation
$\frac{d\mathcal{X}}{dx^1}=\mathcal{A}\mathcal{X}$, and
$\mathcal{X}_0$ be a special solution of the nonhomogeneous equation
$\frac{d\mathcal{X}}{dx^1}=\mathcal{A}\mathcal{X}+\mathcal{F}$, then
to show existence, we need to find $3n$ numbers $c_1,\cdots, c_{3n}$
such that
\begin{eqnarray*}\mathcal{X}=\sum_{j=1}^{3n}c_j\mathcal{X}_j+\mathcal{X}_0
\end{eqnarray*}
satisfies $\mathcal{B}\mathcal{X}=0,$ or equivalently,
\begin{eqnarray}
(\mathcal{B}\mathcal{X}_1,\cdots
,\mathcal{B}\mathcal{X}_{3n})(c_1,\cdots,
c_{3n})^t=-\mathcal{B}\mathcal{X}_0.
\end{eqnarray}
This is a linear algebraic system, and it is well known that it is
uniquely solvable if and only if for $\mathcal{X}_0=0,$ then
$c_1=\cdots=c_{3n}=0.$ However, this is guaranteed by uniqueness we
proved in step 2.1.1.

Note that all the coefficients in the left side of \eqref{316}
belong to  $C^3$ and the right side of \eqref{316} belongs to $C^1$,
so $u^{N,\epsilon}\in C^4.$

{\it Step 2.2.\ Solution of problem \eqref{3009}--\eqref{3011}.} By
the uniform estimate \eqref{320}, there exists a subsequence
$\{u^{N_j,\epsilon}\}$ converges weakly in $H^1$ to a $u^\epsilon\in
H^1$ as $j\rightarrow\infty,$ and $\p_{11}u^{N_j,\epsilon}$
converges weakly in $L^2$ to $\p_{11}u^\epsilon\in L^2$. We show
that $u^\epsilon$ is a weak solution of problem
\eqref{3009}--\eqref{3011}.

For $\chi_j(x^1)\in C_0^{\infty}([-1,1])$, by multiplying it to
\eqref{313}, summing for $j$ from $1$ to $N$, and integrating with
respect to $x^1$ on $[-1,1]$, one has
\[
\int_{M}\{[\p_1(k\chi^N)+\p_2(b\chi^N)+\alpha\chi^N
-\epsilon\p_{11}\chi^N]\p_1u^{N,\epsilon}
+\p_2(a\chi^N)\p_2u^{N,\epsilon}+f_\epsilon\chi^N\}\,dx^1dx^2=0
\]
after integration by parts, where $\chi^N=\sum_{j=1}^N\chi_jY_j.$
Supposing $\chi^N\rightarrow \chi$ strongly in $H^2$, let
$N\rightarrow\infty$ in the above equality, we have
\begin{eqnarray}&&\label{3022}
\int_{M}\{[\p_1(k\chi)+\p_2(b\chi)+\alpha\chi
-\epsilon\p_{11}\chi]\p_1u^{\epsilon}
+\p_2(a\chi)\p_2u^{\epsilon}+f\chi\}\,dx^1dx^2\nonumber\\
&&\quad=0
\end{eqnarray}
for all $\chi\in H^2\cap H^1_0.$ Therefore $u^\epsilon$ is a weak
solution to \eqref{3009}--\eqref{3011}.

Next we show $u^\epsilon$ satisfies \eqref{3010} and \eqref{3011}.
Indeed, since $H^1({M})\hookrightarrow L^2(\p{M})$,
$u^{N,\epsilon}(-1,x^2)=0$ and $u^{N,\epsilon}\rightharpoonup
u^\epsilon$ weakly in $H^1$ indicate that $u^\epsilon(-1,x^2)=0.$
Since $\epsilon$ is fixed presently, and
$W^{1,2}([-1,1];L^2(\mathbf{S}^1))\hookrightarrow
C([-1,1];L^2(\mathbf{S}^1))$ (see Theorem 2, p.286 of \cite{E}),
\eqref{320} implies that $\p_1u^{N,\epsilon}\in
C([-1,1];L^2(\mathbf{S}^1))$. Therefore by \eqref{314} and
\eqref{315}, we obtain that \eqref{3010} and \eqref{3011} hold.

{\it Step 3.\ Existence of weak solution of problem \eqref{201} and
\eqref{202}.}

{\it Step 3.1. Weak solution.} Now for $\epsilon\in (0,1]$ we have
constructed weak solutions of problem \eqref{3009}--\eqref{3011}
which also satisfy \eqref{320}. Then there is a subsequence
$\{u^{\epsilon_j}\}$ converges weakly in $H^1$ to a $u\in H^1$ as
$\epsilon_j\rightarrow0.$ Obviously we have
\begin{eqnarray}\label{8-1-18-1-38}
\norm{u}_1\le C\norm{f}_0.\end{eqnarray}
 We claim $u$ is a weak
solution of problem \eqref{201} and \eqref{202}. To verify this, let
us take $\epsilon=\epsilon_j\rightarrow0$ in \eqref{3022} for any
$\chi\in H^2\cap H^1_0$. Note that
\begin{eqnarray*}
\left|\int_{{M}}\epsilon\p_{11}\chi\p_1u^\epsilon\,dx^1dx^2\right|
\le C \epsilon \norm{f}_0\rightarrow0,
\end{eqnarray*}
we get
\begin{eqnarray}\label{3023}
\int_{M}\{[\p_1(k\chi)+\p_2(b\chi)+\alpha\chi ]\p_1u
+\p_2(a\chi)\p_2u+f\chi\}\,dx^1dx^2=0.
\end{eqnarray}
By approximation this holds for all $\chi\in H^1_0.$

{\it Step 3.2. Boundary conditions.} The fact that $u$ satisfies the
boundary condition $u=0$ on $\Sigma^{-1}$ can be deduced by the same
argument as in step 2.2.

Next, let $M'_\sigma={M}\cap\{(x^1,x^2):1-\sigma<x^1<1\}.$ Since
$k>0$ on $\Sigma^1$ by assumption, for rather small $\sigma$, the
operator $L$ is elliptic in $M'_{3\sigma}.$ We claim that
\begin{eqnarray}\label{h5}
\norm{\p_1u^\epsilon}_{H^1(M'_{2\sigma})}\le C\norm{f}_0.
\end{eqnarray}
If this is true, then clearly $\p_1u=0$ on $\Sigma^1.$

Now we prove \eqref{h5}. Let $\eta\in C^\infty([-1,1])$ such that
\begin{eqnarray}
{\eta(x^1)=}\begin{cases}
0, &-1<x^1\le 1-3\sigma, \label{eq:a1}\\
e^{\mu x^1}, &1-2\sigma\le x^1\le1,\label{eq:a2}
\end{cases}
\end{eqnarray}
where $\mu$ is a large positive constant such that $\p_1(a\eta)\ge0$
on $\Sigma^1.$

Denote $u^{N,\epsilon}$ by $w$. Multiplying \eqref{313} by
$\eta\p_{11}(X_j^{N,\epsilon})$, summing up for $j$ from $1$ to $N$,
and integrating the equality on $[-1,1]$, we have
\begin{eqnarray}
&&\int_{{M}}f\eta \p_{11}w\, dx^1dx^2\nonumber\\
&=&\int_{{M}} \big\{\p_1\big[\frac12 \epsilon\eta(\p_{11}w)^2+\frac
12\p_1(\eta
a)(\p_2w)^2-\frac12\alpha\eta(\p_1w)^2-a\eta\p_2w\p_{12}w\big]\nonumber\\
&&+\p_2(a\eta\p_2w\p_{11}w)+(k\eta-\frac12\epsilon\p_1\eta)(\p_{11}w)^2+b\eta\p_{11}w
\p_{12}w+\frac12\p_1(\alpha\eta)(\p_1w)^2\nonumber\\
&&-\p_2(a\eta)\p_2w\p_{11}w-\frac12\p_{11}(a\eta)(\p_2w)^2
+a\eta(\p_{12}w)^2\big\}\,dx^1dx^2\nonumber\\
&\ge&\int_{{M}}\eta[(k-\nu-\varrho)(\p_{11}w)^2+(a-\nu)(\p_{12}w)^2]\,dx^1dx^2-
C_{\tau,\varrho}\norm{w}_1^2.\label{h7}
\end{eqnarray}
Here we used the fact the $\p_1w=\p_{12}w=0$ on $\Sigma^1,$ and
\eqref{320} to control the term
$-(1/2)\epsilon\p_1\eta(\p_{11}w)^2$. Hence by applying Young
inequality to the left side of the above inequality and choosing
$\nu,\ \varrho$ small enough, we have
\begin{eqnarray}\label{h6}
\int_{{M}}\eta[(\p_{11}w)^2+(\p_{12}w)^2]\,dx^1dx^2\le
C\norm{f}_0^2.
\end{eqnarray}
This implies \eqref{h5} by taking $N\rightarrow\infty$.

{\it Step 4.\ $H^2$ Regularity.}

{\it Step 4.1.\ Regularity in elliptic region.} Let
$M_\sigma={M}\cap\{(x^1,x^2):-1<x^1<-1+\sigma\}.$ Since $k>0$ on
$\Sigma^{-1},$ we may choose $\sigma>0$ small such that the operator
$L$ is elliptic in $M_{2\sigma}.$ By standard theory of elliptic
equations (i.e., the interior estimate and the estimate near
boundary of elliptic equations) we have
\begin{eqnarray}\label{3024}
&&\norm{u}_{H^2(M_{2\sigma})}\le C(\norm{u}_0+\norm{f}_0)\le
C\norm{f}_0,\\
&&\norm{u}_{H^2(M'_{2\sigma})}\le C\norm{f}_0.\label{3025}
\end{eqnarray}

Next we derive an estimate of $w=u^{N,\epsilon}$ similar to
\eqref{h6} near $\Sigma^{-1}.$ Let $\xi=\xi(x^1)\in
C^{\infty}_0([-1,1])$ satisfy
\begin{eqnarray} \label{eq:xi1}
\xi(x^1)=\left\{ \begin{aligned}
         0, &\ x^1\in[-1,-1+\sigma/4]\cup[-1+2\sigma, 1],\\
                  1, &\  x^1\in[-1+\sigma/2,-1+3\sigma/2].
                          \end{aligned} \right.
                          \end{eqnarray}
By multiplying $\xi\p_{11}(X_j^{N,\epsilon})$ to \eqref{313},
summing up for $j$ from $1$ to $N$ and then integrating on $[-1,1]$
with respect to $x^1$, by the first equality in \eqref{h7} we have
similarly
\begin{eqnarray}
&&\int_{{M}}f\xi \p_{11}w\, dx^1dx^2\nonumber\\
&\ge&(\delta/16)\int_{{M}}\xi[(\p_{11}w)^2+(\p_{12}w)^2]\,dx^1dx^2-
C\norm{w}_1^2,\nonumber
\end{eqnarray}
hence
\begin{eqnarray}\label{h8}
\int_{{M}}\xi[(\p_{11}w)^2+(\p_{12}w)^2]\,dx^1dx^2\le C\norm{f}_0^2.
\end{eqnarray}

{\it Step 4.2.\ Regularity in mixed--type region.}  Let $\vartheta
\in C_0^\infty([-1,1])$ such that
\begin{eqnarray*}
0\leq\vartheta \leq 1\quad \hbox{and} \quad \vartheta
(x^1)=\begin{cases}
0, & -1\leq  x^1\leq-1+\sigma/2,\\
1, & -1+3\sigma/2\leq x^1\leq 1-3\sigma/2,\\
0, & 1-\sigma/2\leq x^1\leq 1.
\end{cases}
\end{eqnarray*}
Multiplying \eqref{313} by $\vartheta \p_{111}(X_j^{N,\epsilon})$,
summing up for $j$ from $1$ to $N$, and integrating the equality on
$[-1,1]$, we get
\begin{eqnarray}\label{hh}
&&-\int_{{M}}\p_1(f\vartheta )\p_{11}w\,dx^1dx^2=\int_{{M}}
\{(\alpha-\frac12\p_1k)\vartheta (\p_{11}w)^2-\frac12k\vartheta '(\p_{11}w)^2\nonumber\\
&&\qquad-\frac32\p_1a\vartheta (\p_{12}w)^2-\frac32 a\vartheta
'(\p_{12}w)^2 -\p_{11}(a\vartheta )\p_2w\p_{12}w
-\frac12\p_{11}(\alpha\vartheta )(\p_1w)^2\nonumber\\
&&\qquad\ \ +\p_{12}(a\vartheta )\p_2w\p_{11}w-\p_1(b\vartheta
)\p_{12}w\p_{11}w
+\epsilon\vartheta (\p_{111}w)^2\nonumber\\
&&\qquad\quad+\frac12 \p_2(b\vartheta )(\p_{11}w)^2+\p_2(a\vartheta
)\p_{12}w\p_{11}w
 \}\,dx^1dx^2.
\end{eqnarray}
Here we may estimate those terms involving $\vartheta ',\,\vartheta
''$ by the estimates \eqref{h6} and \eqref{h8}, since
$(\supp\,\,\vartheta')\subset [-1+\sigma/2,-1+3\sigma/2]\cup
[1-3\sigma/2,1-\sigma/2]$. Therefore
\begin{eqnarray}
-\int_M [\frac12k\vartheta'(\p_{11}w)^2+\frac32
a\vartheta'(\p_{12}w)^2]\,dx^1dx^2 \leq C\|f\|_0.
\end{eqnarray}
Hence by H\"older inequality and Young inequality, the right side of
\eqref{hh} is bounded below by
\begin{eqnarray}
(\delta/4)\int_{{M}}\vartheta
|D(\p_1w)|^2\,dx^1dx^2-C(\norm{w}_1^2+\|f\|_1^2).
\end{eqnarray}
Therefore, it is easy to get
\begin{eqnarray*}
\int_{{M}}\vartheta |D(\p_1w)|^2\,dx^1dx^2\le C\norm{f}_1^2.
\end{eqnarray*}
with $C$ independent of $N, \epsilon$. Letting $N\rightarrow \infty$
and $\epsilon\rightarrow 0$, we obtain that
\begin{eqnarray*}
\int_{{M}}\vartheta |D(\p_1u)|^2\,dx^1dx^2\le C\norm{f}_1^2.
\end{eqnarray*}
By \eqref{8-1-18-1-38}, \eqref{3024}, \eqref{3025} and the above
inequality, we obtain that
\begin{eqnarray*}
\|\p_1u\|_1\le C\norm{f}_1.
\end{eqnarray*}

Since $a\ge\delta>0$ in $M$, by \eqref{h1} we may estimate
$\p_{22}u$. Hence we conclude that $\norm{u}_2\le C\norm{f}_1$.

\begin{rem}\label{rem301}
We observe that in deriving the $H^2$ estimate, we used the
assumption that $D^2a, D^2\alpha\in L^\infty$, but just required
that $Dk, Db\in L^\infty.$
\end{rem}

{\it Step 5. Higher regularity.} The regularity in \eqref{204} for
$s=3, 4$ in the elliptic region is obvious (also can be obtained by
multiplier $\xi \p_{1}^{k}w$,$\eta \p_{1}^{k}w$  for integer $k=4,
6$), and in the mixed-type region can be obtained in the same
fashion as those to \eqref{h6} by multiplier $\vartheta \p_{1}^{k}w$
for integer $k=5, 7$.

For example, for $H^3$ estimate, by the above energy estimate
technique, we have
\begin{eqnarray}
&&-\int_{M}\p_{11}(f\vartheta)\p_{111}w\,dx^1dx^2=\int_{M}\epsilon\vartheta(\p_{111}w)^2\,dx^1dx^2\nonumber\\
&&\quad+\int_{M}\{[\alpha\vartheta-\frac32\p_1(k\vartheta)+\underline{\frac12\p_2(b\vartheta)}](\p_{111}w)^2-\frac52\p_1
(a\vartheta)(\p_{112}w)^2\}\,dx^1dx^2\nonumber\\
&&\qquad+\int_{M}[\p_2(a\vartheta)\p_{112}w\p_{111}w-\underline{2\p_1(b\vartheta)\p_{112}w\p_{111}w}]\,dx^1dx^2\nonumber\\
&&\qquad+\int_{M}[-\p_{11}(a\vartheta)\p_{22}w\p_{111}w+2\p_{12}(a\vartheta)\p_{12}w\p_{111}w-2
\p_{11}(a\vartheta)\p_{12}w\p_{112}w\nonumber\\
&&\qquad\ \
+\p_{11}(\alpha\vartheta)\p_1w\p_{111}w-\p_{11}(k\vartheta)\p_{11}w\p_{111}w-\p_{11}(b\vartheta)
\p_{12}w\p_{111}w]\,dx^1dx^2\nonumber\\
&&\qquad\
\ -\int_{M}\p_{11}(\alpha\vartheta)(\p_{11}w)^2\,dx^1dx^2\nonumber\\
&\ge&\int_{M}\vartheta\{(\delta-\nu-\sigma)[(\p_{111}w)^2+(\p_{112}w)^2]\}\,dx^1dx^2-C_\sigma\norm{f}_2^2.
\end{eqnarray}
An important fact is that the terms with underlines involve only
$Db,$ so we may use \eqref{2005} (and $H^3\hookrightarrow C^1$) to
get the inequality. Therefore by choosing $\sigma$ small we get
\begin{eqnarray}
\int_{M}\vartheta[(\p_{111}w)^2+(\p_{112}w)^2]\,dx^1dx^2\le
C\norm{f}_2^2.
\end{eqnarray}
Combing an estimate obtained similarly in the elliptic domain, we
have
\begin{eqnarray}
\norm{\p_{11}w}_1\le C\norm{f}_2.
\end{eqnarray}
Letting $N\rightarrow\infty$ and $\epsilon\rightarrow0,$ the weak
limit $u$ of the sequence $\{w=u^{N,\epsilon}\}$, which is a
solution of \eqref{h1}, also satisfies
\begin{eqnarray}
\norm{\p_{11}u}_1\le C\norm{f}_2.
\end{eqnarray}
By \eqref{h1}, we have
\begin{eqnarray}
\norm{\p_{22}u}_1&\le&
C\norm{f}_1+C\norm{u}_2+C\norm{\p_{11}u}_1+\norm{\frac
ba\p_{12}u}_1\nonumber\\
&\le&C\norm{f}_2+\norm{\frac ba\p_{112}u}_0+\norm{\frac
ba\p_{122}u}_0
\end{eqnarray}
for a constant $C$ depending on $\delta.$ By \eqref{2005}, we have
$\norm{\frac
ba\p_{122}u}_0\le(\nu/\delta)\norm{\p_{22}u}_1<\norm{\p_{22}u}_1.$
Also, there holds $\norm{\frac ba\p_{112}u}_0\le
C\norm{\p_{11}u}_1$. So we get
\begin{eqnarray}
\norm{\p_{22}u}_1\le C\norm{f}_2
\end{eqnarray}
and the $H^3$ regularity
\begin{eqnarray}
\norm{u}_3\le \norm{u}_2+\norm{\p_{11}u}_1+\norm{\p_{22}u}_1\le
C\norm{f}_2.
\end{eqnarray}

Similar argument works for $H^4$ estimate. This finishes the proof.
\hfill$\square$

\begin{rem}\label{rem302}
We see here that we can not apply directly Theorem \ref{thm201} to
\eqref{105} since $\p_1a\le-\delta<0$ in \eqref{h0} does not hold in
${M}$. However, we need this assumption to control the mixed
derivative term $b\p_{12}u$.

An important observation to Theorem \ref{thm201} is that its
conditions are not invariant under multiplication of a positive
function to the mixed type equation \eqref{201}! Therefore it is
expected to find an appropriate multiplier to \eqref{105} such that
Theorem \ref{thm201} works.
\end{rem}

Now we  choose $h(x^1)=e^{-\mu{x^1}}$ with $\mu>0$ a small constant
(depending only on $\delta_2$ and $\norm{k}_{L^{\infty}}$), which is
a bounded smooth positive function in $\mathcal{M}$, then it is easy
to check that there is a positive constant $\delta^*$ such that
\begin{eqnarray}
&h'(x^1)\le-\delta^*<0 & \text{in}\ \ {M},\label{20501}\\
&2h(x^1)\alpha(x^1)+\p_1(h(x^1)k(D\varphi_b))\ge\delta^*>0
 & \text{in}\ \ {M},\label{206}\\
 &2\alpha h-l\p_1(h k(D\varphi_b))\ge\delta^*>0 &\text{in} \ \
{M},\label{211}
\end{eqnarray}
where $l>0$ is less than a fixed number (say, $l\le6$).

\section{Solvability of nonlinear problem and main result}

Now we prove the following stability result of transonic potential
flows in $\mathcal{M}$.
\begin{thm}\label{thm301}
Let $\varphi_b\in C^5$ be a background transonic flow in
$\mathcal{M}$. Then there exist  positive constant $C$ and
$\varepsilon_0$ depending only on $\varphi_b$ such that if
\begin{eqnarray}\label{301}
\varphi=\varphi_b+g(x^2) \qquad \text{on}\ \ \Sigma^{-1}
\end{eqnarray}
for any $g\in H^5(\mathbf{S}^1)$ and $\norm{g}_{5}\le
\varepsilon\le\varepsilon_0$, then problem \eqref{101}, \eqref{102}
and \eqref{301} has uniquely one solution $\varphi$ with
\begin{eqnarray}\label{302}
\norm{\varphi-\varphi_b}_4\le C\varepsilon.
\end{eqnarray}
\end{thm}

This result shows that the subsonic--supersonic flow is stable under
perturbations of the velocity potential function at the entry. More
physically, since $\p_2\varphi/n(-1)$ is the velocity of the flow
along $x^2$ direction, so $g(x^2)$ measures the perturbation of the
tangential velocity along the entry. So we may claim that the
special transonic flow is stable under small variation of the
tangential velocity at the entrance.

\pf The proof is based on a nonlinear iteration scheme.

{\it Step 1.} We define the iteration set as
\begin{eqnarray*}
E_\kappa=\{\varphi\in H^4({M}):
\norm{\varphi-\varphi_b}_4\le\kappa\le\kappa_0\},
\end{eqnarray*}
where $\kappa_0$ is a small positive constant to be specified later.
It is straightforward to check that there hold
\begin{eqnarray}
&&\norm{f(D\varphi)}_3\le C_0\kappa^2, \label{30401}\\
&&\norm{f(D\varphi^{(1)})-f(D\varphi^{(2)})}_2\le
C_0\kappa\norm{\varphi^{(1)}-\varphi^{(2)}}_3 \label{30501}
\end{eqnarray}
for any $\varphi, \varphi^{(1)}, \varphi^{(2)}\in E_\kappa.$

{\it Step 2.} Let
\begin{eqnarray}\label{304}
(h\cdot
M_\varphi)(\hat{\varphi})&:=&h(x^1)k(D\varphi)\p_{11}\hat{\varphi}
+h(x^1)b(D\varphi)\p_{12}\hat{\varphi}\nonumber\\
&&\qquad
+h(x^1)\p_{22}\hat{\varphi}-h(x^1)\alpha(x^1)\p_1\hat{\varphi}\nonumber\\
&=&h(x^1)f(D{\varphi}),
\end{eqnarray}
where $\hat{\varphi}=\varphi-\varphi_b,$ and $f(D{\varphi})$ is
defined by \eqref{110}. By considering $\phi=\hat{\varphi}-g(x^2)$
as the unknown, problem \eqref{101}\eqref{102} and \eqref{301} is
equivalent to the following problem:
\begin{eqnarray}
&(hM_\varphi)(\phi)=h(x^1)f(D{\varphi})+(hM_\varphi)(g(x^2))
&\text{in}\ \ {M},\\
&\phi=0 &\text{on}\ \ \Sigma^{-1}.
\end{eqnarray}
Here $\varphi=\phi+g(x^2)+\varphi_b,\
\hat{\varphi}=\varphi-\varphi_b$ and note that
\begin{eqnarray}
&&\norm{M_\varphi(g)}_3\le C_0\varepsilon,\\
&&\norm{M_{\varphi^{(1)}}(g)-M_{\varphi^{(2)}}(g)}_2\le
C_0\varepsilon\norm{\varphi^{(1)}-\varphi^{(2)}}_3
 \end{eqnarray} for
$\varphi, \varphi^{(1)}, \varphi^{(2)}\in E_\kappa.$

{\it Step 3.} Now by \eqref{20501}--\eqref{211}, we may choose
$\kappa_0$ so small such that for any $\varphi\in E_\kappa$, there
hold in ${M}$ the following inequalities:
\begin{eqnarray}
&h(x^1)\ge\delta^*/2>0,\label{307}\\
&2h(x^1)\alpha(x^1)+\p_1(h(x^1)k(D\varphi))\ge\delta^*/2>0,\label{308}\\
&\p_1h(x^1)\le-\delta^*/2<0,\label{309}\\
&2\alpha h-l\p_1(h k(D\varphi))\ge\delta^*/2>0&
l=0,1,\cdots,5,\label{310}\\
&\norm{b(D\varphi)}_{3}\le \nu^*=\delta^*/4.\label{b1}
\end{eqnarray}

{\it Step 4.} Then for any $\varphi\in E_\kappa,$ we solve the
following linear problem of $\bar{\phi}$ :
\begin{eqnarray}
&(hM_\varphi)(\bar{\phi})=h(x^1)f(D{\varphi})+(hM_\varphi)(g(x^2))
&\text{in}\ \ {M},\label{416}\\
&\bar{\phi}=0 &\text{on}\ \ \Sigma^{-1}.\label{417}
\end{eqnarray}
By \eqref{307}--\eqref{b1} and Theorem \ref{thm201}, and the
analysis in step 6 below, there exists uniquely one solution
$\bar{\phi}\in H^1$ and it satisfies
\begin{eqnarray*}
\norm{\bar{\phi}}_4\le C_0(\kappa^2+\varepsilon).
\end{eqnarray*}
Now choosing $\varepsilon_0\le1/(8C_0^2)$ and
$\kappa=4C_0\varepsilon\le\kappa_0$ (that is, $C=4C_0$), we get a
$\bar{\varphi}=\bar{\phi}+g+\varphi_b$ with
$\norm{\bar{\varphi}}_4\le\kappa.$ Therefore we established a
mapping $T: \varphi\mapsto\bar{\varphi}$ on $E_\kappa.$

{\it Step 5.} Next we will show that $T$ is contractive on
$E_\kappa$ in $H^3$ norm.

Let $\varphi^{(i)}\in E_\kappa,$
$T(\varphi^{(i)})=\bar{\varphi}^{(i)},$ and
$\bar{\phi}^{(i)}=\varphi^{(i)}-g-\varphi_b,\,\, i=1,2.$ Then
$\bar{\phi}^{(1)}-\bar{\phi}^{(2)}$ satisfies the following problem
\[
(hM_{\varphi^{(1)}})(\bar{\phi}^{(1)}-\bar{\phi}^{(2)})
=-\big[hM_{\varphi^{(1)}}(\phi^{(2)})-hM_{\varphi^{(2)}}(\phi^{(2)})\big]
+h(x^1)\big[f(D\varphi^{(1)})-f(D\varphi^{(2)})\big]
\]
\[
\qquad\qquad+(hM_{\varphi^{(1)}})(g)-(hM_{\varphi^{(2)}})(g)\qquad\qquad\
\ \ \ \text{in}\ \ \ M,
\]
\[
\bar{\phi}^{(1)}-\bar{\phi}^{(2)}=0\qquad\qquad\qquad\qquad\qquad\qquad\qquad\qquad
\text{on}\ \
 \Sigma^{-1}.
\]
Note that
\begin{eqnarray}
\norm{M_{\varphi^{(1)}}\phi^{(2)}-M_{\varphi^{(2)}}\phi^{(2)}}_2\le
C_0 C\varepsilon\norm{\varphi^{(1)}-\varphi^{(2)}}_3,
\end{eqnarray}
then Theorem \ref{thm201} and the results in step 6 below yield
\begin{eqnarray}
\norm{\bar{\varphi}^{(1)}-\bar{\varphi}^{(2)}}_3\le
C_1\varepsilon\norm{\varphi^{(1)}-\varphi^{(2)}}_3.
\end{eqnarray}
By choosing $\varepsilon_0$ further small, we obtain contraction.
Then by a simple generalized Banach fixed point theorem, we proved
Theorem \ref{thm301}.

{\it Step 6. Solving $H^4$ solution of \eqref{416} and \eqref{417}
with lower regular coefficients.}

 For simplicity, we may write
\begin{eqnarray}\label{423}
M_\varphi(\bar{\phi}):=[\bar{k}+O_1(D\hat{\varphi})]\p_{11}\bar{\phi}
+O_2(D\hat{\varphi})\p_{12}\bar{\phi}
+\p_{22}\bar{\phi}-{\alpha}\p_1{\bar{\phi}}=F(D{\varphi}),
\end{eqnarray}
where $\bar{k}=k(D\varphi_b),\
O_1(D\hat{\varphi})=k(D\varphi)-k(D\varphi_b),\
O_2(D\hat{\varphi})=b(D\varphi)$, $F=f(D\varphi)+M_\varphi g$, and
there holds
\begin{eqnarray}\label{424}
\norm{O_i}_3\le C\kappa\le C\kappa_0
\end{eqnarray}
for $i=1,2$ and $\varphi\in E_\kappa.$ Since $h$ is smooth and
bounded away from zero, the solvability of \eqref{423}\eqref{417} is
equivalent to \eqref{416}\eqref{417}.

By Sobolev embedding theorem  we have $H^3\hookrightarrow C^1$, so
due to  Remark \ref{rem301} in Theorem \ref{thm201}, \eqref{423},
\eqref{417} has uniquely one solution $\bar{\phi}$ and
\begin{eqnarray}\label{425}
\norm{\bar{\phi}}_2\le C\norm{F}_1.
\end{eqnarray}
However, we can not infer from Theorem \ref{thm201} directly the
existence of $H^4$ solutions and estimate like
$\norm{\bar{\phi}}_4\le C\norm{F}_3$. So we need the following
apriori estimates and approximation argument.

{\it Step 6.1. An apriori $H^4$ estimate.} Now suppose there is a
solution $\bar{\phi}\in H^4$ to problem \eqref{423}\eqref{417}, we
show that if $\kappa_0$ is sufficiently small, then there holds the
following estimate
\begin{eqnarray}\label{426}
\norm{\bar{\phi}}_p\le C\norm{F}_{p-1},\quad \text{for}\,\,\,p=3,4.
\end{eqnarray}

The proof of the case $p=3$ is similar to the following $p=4$ case
and even more simple (which should use the fact \eqref{425}). Now we
suppose \eqref{426} holds for $p=3$ and to prove $p=4$ case.

First, for small $\sigma$ and $\kappa_0,$ the operator $M_\varphi$
is elliptic in $M_{3\sigma}.$ By the regularity theory of elliptic
equation, we have
\begin{eqnarray}\label{427}
\norm{\bar{\phi}}_{H^4(M_{3\sigma})}\le C\norm{F}_3.
\end{eqnarray}

Now let $\zeta \in C^\infty([-1,1])$  satisfy
\begin{eqnarray*}
0\leq\zeta \leq 1\quad \hbox{and} \quad \zeta (x^1)=\begin{cases}
0, & -1\leq x^1\leq -1+\sigma,\\
1, & -1+2\sigma\le x^1\le1.
\end{cases}
\end{eqnarray*}
Denote $\bar{u}=\zeta\bar{\phi}.$ Then it satisfies the following
equation
\begin{eqnarray}\label{428}
M_\varphi(\bar{u})=\bar{F}:=\zeta
F+\zeta''k\bar{\phi}+\zeta'(2k\p_1\bar{\phi}+b\p_2\bar{\phi}-\alpha\bar{\phi}).
\end{eqnarray}
Since $\supp\, \zeta'\subset(-1+\sigma, -1+2\sigma)$, by
\eqref{427}, we have
\begin{eqnarray}\label{429}
\norm{\bar{F}}_3\le C\norm{F}_3.
\end{eqnarray}

Then by differentiating \eqref{428} with respect to $x^1$ twice, we
get $w=\p_{11}\bar{u}$ satisfies
\begin{eqnarray}\label{430}
M_\varphi(w)+2\p_1\bar{k}\p_1w
&=&\p_{11}\bar{F}+[-\p_{11}\bar{k}\p_{11}\bar{u}+\p_{11}\alpha\p_1\bar{u}
+2\p_1\alpha\p_{11}\bar{u}]\nonumber\\
&&-[\p_{11}O_1\p_{11}\bar{u}+\p_{11}O_2\p_{12}\bar{u}
+2\p_1O_1\p_{111}\bar{u}+2\p_1O_2\p_{112}\bar{u}]\nonumber\\
&:=&\tilde{F}
\end{eqnarray} as well as $w=0$ on $\Sigma^{-1}$ by
the cut-off function $\zeta.$

We note that the operator $M_\varphi+2\p_1\bar{k}\p_1$ in the left
hand side of the above equation also satisfies the assumptions of
Theorem \ref{thm201} (by multiplying a suitable positive function).
So we have
\begin{eqnarray}
\norm{w}_2\le C \|\tilde{F}\|_1.
\end{eqnarray}
We now estimate $\tilde{F}$ term by term.

(1) Obviously by \eqref{429} we have
\begin{eqnarray}
\norm{\p_{11}\bar{F}}_1\le \norm{\bar{F}}_3\le C\norm{F}_3.
\end{eqnarray}

(2) Since $\bar{k}, \alpha\in C^4$, we get
\begin{eqnarray}
\norm{-\p_{11}\bar{k}\p_{11}\bar{u}+\p_{11}\alpha\p_1\bar{u}
+2\p_1\alpha\p_{11}\bar{u}}_1\le C\norm{\bar{u}}_3\le C\norm{F}_2.
\end{eqnarray}
Here we use \eqref{425} and \eqref{426} ($p=3$) for the second
inequality.

(3) We first recall the inequality
\begin{eqnarray}
\norm{uv}_1\le C\norm{u}_2\norm{v}_1
\end{eqnarray}
provided $u\in H^2,\ v\in H^1$ (see \cite{K2}, p.73). Then by
\eqref{424}
\begin{eqnarray*}
&&\norm{\p_{11}O_1\p_{11}\bar{u}}_1\le
C\norm{O_1}_3\norm{\bar{u}}_4\le C\kappa\norm{\bar{u}}_4,\\
&&\norm{\p_{11}O_2\p_{12}\bar{u}}_1\le
C\norm{O_2}_3\norm{\bar{u}}_4\le C\kappa\norm{\bar{u}}_4,\\
&&\norm{\p_1O_1\p_{111}\bar{u}}_1\le
C\norm{O_1}_3\norm{\bar{u}}_4\le
C\kappa\norm{\bar{u}}_4,\\
&&\norm{\p_1O_2\p_{112}\bar{u}}_1\le
C\norm{O_2}_3\norm{\bar{u}}_4\le C\kappa\norm{\bar{u}}_4.
\end{eqnarray*}

In all, we get $\|\tilde{F}\|_1\le
C(\norm{F}_3+\kappa\norm{\bar{u}}_4),$ hence
\begin{eqnarray}\label{439}
\norm{\p_{11}\bar{u}}_2\le C(\norm{F}_3+\kappa\norm{\bar{u}}_4).
\end{eqnarray}
By \eqref{428}, we may then estimate
\begin{eqnarray*}
\norm{\p_{22}\bar{u}}_2\le
C(\norm{F}_3+\kappa\norm{\p_{12}\bar{u}}_2+\kappa\norm{\bar{u}}_4)\le
C(\norm{F}_3+\kappa\norm{\bar{u}}_4).
\end{eqnarray*}
So there holds
\begin{eqnarray*}
\norm{\bar{u}}_4\le
C(\norm{\p_{11}\bar{u}}_2+\norm{\p_{22}\bar{u}})+\norm{\bar{u}}_3\leq
C(\norm{F}_3+\kappa\norm{\bar{u}}_4).
\end{eqnarray*}
By choosing $\kappa_0$ small, we can deduce that
$\norm{\bar{u}}_4\leq C\norm{F}_3$. Combing this with \eqref{427},
we can get \eqref{426} for the case $p=4$.

{\it Step 6.2. Existence of  $H^4$ solution by approximation.}

For $\varphi\in E_\kappa,$ $O_i\in H^3, i=1, 2,$ we approximate
$O_i$ by $\{O_i^{(l)}\}_{l=1}^\infty\subset C^{4}$ such that
$O_i^{(l)}\rightarrow O_i(D\hat{\varphi})$ strongly in  $H^3$ (so
\eqref{424} holds). By Theorem \ref{201}, the problem
\begin{eqnarray*}
&&M^{(l)}_\varphi(\bar{\phi}):=[\bar{k}+O_1^{(l)}]\p_{11}\bar{\phi}
+O_2^{(l)}\p_{12}\bar{\phi}
+\p_{22}\bar{\phi}-{\alpha}\p_1{\bar{\phi}}=F(D{\varphi}),\\
&&\bar{\phi}=0 \quad \text{on}\quad \Sigma^{-1}
\end{eqnarray*}
has uniquely one solution $\bar{\phi}^{(l)}\in H^4.$ By the apriori
estimate  \eqref{426} we have $\norm{\bar{\phi}^{(l)}}_4\le
C\norm{F}_3$ for $C$ independent of $l.$ So there is a $H^4$ weak
limit $\bar{\phi}\in H^4$ of this sequence of approximate solutions.
Then clearly $\bar{\phi}$ is a $H^4$ solution of \eqref{423} and
\eqref{417}, and by Theorem \ref{thm201}, this solution is unique.

The proof of Theorem \ref{thm301} is then finished.


\end{document}